\documentclass[10pt]{article}
\usepackage{latexsym,amssymb,amsmath,amsthm,amsfonts}
\usepackage{graphicx}
\usepackage{enumerate}
\usepackage{color}

\theoremstyle{plain}
\newtheorem{thm}{Theorem}[section]
\newtheorem{prop}[thm]{Proposition}
\newtheorem{lem}[thm]{Lemma}

\theoremstyle{definition}
\newtheorem{definition}[thm]{Definition}

\numberwithin{equation}{section}

\newcommand{\xy}{\langle x,y\rangle}

\newcommand{\pppp}[4]%
  {\frac{\partial^3{#1}}{\partial{#2}\partial{#3}\partial{#4}}}

\newcommand{\p}{\phi}

\newcommand{\gab}{\alpha\phi \Big(b^2,\frac{\beta}{\alpha}\Big )}
\newcommand{\pt}{\phi_2}
\newcommand{\po}{\phi_1}
\newcommand{\ptt}{\phi_{22}}
\newcommand{\pot}{\phi_{12}}

\renewcommand{\a}{\alpha}
\renewcommand{\b}{\beta}
\newcommand{\ab}{(\alpha,\beta)}

\begin{document}
\title{On general $\ab$-metrics with isotropic Berwald curvature}
\footnotetext{\emph{Keywords}:  Finsler metric, general $\ab$-metric, Douglas metric.
\\
\emph{Mathematics Subject Classification}: 53B40, 53C60.}
\author{Hongmei Zhu\footnote{supported in  a doctoral scientific research foundation of Henan Normal University (No.5101019170130)}}
\maketitle


\begin{abstract}
In this paper, we study a class of Finsler metrics called general $\ab$-metrics, which are defined by a Riemannian metric $\a$ and a $1$-form $\b$.
We classify this class of Finsler metrics with isotropic Berwald curvature under certain condition.
\end{abstract}

\section{Introduction}
In Finsler geometry, the Berwald curvature is an important  non-Riemannian quantity. A Finsler metric $F$ on a manifold $M$ is said to be of {\it isotropic Berwald curvature} if its Berwald curvature $B_{j}{}^{i}{}_{kl}$ satisfies
\begin{eqnarray}\label{isotropic}
B_{j}{}^{i}{}_{kl}=\tau(x)(F_{y^{j}y^{k}}\delta^{i}{}_{l}+F_{y^{j}y^{l}}\delta^{i}{}_{k}+F_{y^{l}y^{k}}\delta^{i}{}_{j}+F_{y^{j}y^{k}y^{l}}y^{i}),
 \end{eqnarray}
where $\tau(x)$ is a scalar function on $M$. A Finsler metric is called a {\it Berwald metric} if $\tau(x)=0$. Berwald metrics are just a bit more general than Riemannian and locally
Minkowskian metrics. A Berwald space is that all tangent spaces are linearly isometric to a common Minkowski space.
 
Chen-Shen showed that a Finsler metric $F$ is of isotropic Berwald curvature if and only if it is a Douglas metric with isotropic mean Berwald curvature \cite{Chen-Shen}. Tayebi-Rafie's result tells us that every isotropic Berwald metric is of isotropic $S$-curvature \cite{Tayebi-Rafie}. In \cite{Tayebi-Najafi}, Tayebi-Najafi proved that
isotropic Berwald metrics of scalar flag curvature are of Randers type. Recently, Guo-Liu-Mo have shown that every spherically symmetric Finsler metric of isotropic Berwald curvature is a Randers metric \cite{Guo-Liu-Mo}.

Randers metrics , which are introduced by a physicist G. Randers in 1941 when he studied general relativity,  are an important class of Finsler metrics. It is of the form $F=\a+\b$, where $\a$ is a Riemannian metric and $\b$ is a 1-form. Bao-Robles-Shen show that Randers metric $F=\a+\b$,
where
\begin{eqnarray*}
\a:&=& \frac{\sqrt{\varepsilon^{2}(xu+yv+zw)^{2})+(u^{2}+v^{2}+w^{2})[1-\varepsilon^{2}(x^{2}+y^{2}+z^{2})]}}{1-\varepsilon^{2}(x^{2}+y^{2}+z^{2})},\\
\b:&=& \frac{-\varepsilon(xu+yv+zw)}{1-\varepsilon^{2}(x^{2}+y^{2}+z^{2})}.
\end{eqnarray*}
is of constant $S$-curvature and satisfies $d\beta=0$. Hence, $F$ is of isotropic Berwald curvature\cite{Bao-Robles-Shen, Guo-Liu-Mo}.

Many famous Finsler metrics can be expressed  in the following form
\begin{equation}\label{generalab}
F=\gab,
\end{equation}
where $\a$ is a Riemannian metric, $\b$ is a $1$-form, $b: =\|\beta_x\|_{\alpha}$ and $\p(b^2,s)$ is a smooth function. Finsler metrics in this form are called general $\ab$-metrics \cite{yct-zhm-onan, yu-zhu}. If $\phi=\phi(s)$ is independent of $b^2$, then $F=\a\phi(\frac{\b}{\a})$ is an $\ab$-metric. If $\a=|y|$, $\b=\xy$, then $F=|y|\phi(|x|^2,\frac{\xy}{|y|})$ is the so-called spherically symmetric Finsler metrics \cite{Mo-Zhu, Zhm}. Moreover, general $\ab$-metrics include part of Bryant's metrics \cite{Br2, yct-zhm-onan} and part of generalized fourth root metrics \cite{Li-Shen}.
 In fact, Randers metrics $F=\a+\b$ can also be expressed in the following navigation form
\begin{equation*}
F = \frac{\sqrt{(1-b^2)\bar{\alpha}^2+\bar{\beta}^2 }}{1-\bar{b}^2} +\frac{\bar{\beta}}{1-\bar{b}^2},\label{Randersmetric}
\end{equation*}
where $\bar{\alpha}$ is also a Riemannian metric, $\bar{\b}$ is a $1$-form and $\bar{b}:=\|\bar{\b}_{x}\|_{\bar{\alpha}}$. $(\bar{\a},\bar{\b})$ is called {\it the navigation data } of the Randers metric $F$.
In \cite{Tayebi-Rafie}, Tayebi-Rafie showed that if a Randers metric $F=\a+\b$ is an non-trivial isotropic Berwald metric, then $\bar{\b}$ is a conformal $1$-form with respect to $\bar{\a}$, namely, $\bar{\b}$ satisfies
$$\bar{b}_{i|j}+\bar{b}_{j|i}=c(x)\bar{a}_{ij},$$
where $c(x)\neq0$ and $\bar{b}_{i|j}$ is the covariant derivation of $\bar{\b}$ with respect to $\bar{\a}$.

In this paper, we mainly study general $\ab$-metrics $F=\a \phi(b^{2},s)$ with isotropic Berwald curvature and  show the following classification theorem

\begin{thm}\label{main2} Let $F=\gab$ be a general $\ab$-metric on an $n$-dimensional manifold $M$. Suppose that $\b$ satisfies
\begin{equation}
b_{i|j} = c a_{ij}, \label{bcij}
\end{equation}
where $c=c(x)\neq 0$ is a scalar function on $M$ . If $F$ is of isotropic Berwald curvature, then one of the following holds

(1) $F$ is Riemannian;

(2) $F$ is a Randers metric;

(3) $F$ is a Kropina metric;

 (3) $F$ is a Berwald metric which can be expressed by
\begin{eqnarray*}
 F=\alpha \varphi(\frac{s^{2}}{e^{\int(\frac{1}{b^{2}}-b^{2}t_{2}))db^{2}}+s^{2}\int t_{2}e^{\int(\frac{1}{b^{2}}-b^{2}t_{2})db^{2}}db^{2}})e^{\int(\frac{1}{2}b^{2}t_{2}-\frac{1}{b^{2}})db^{2}}s.
 \end{eqnarray*}
 where $\varphi(\cdot)$ is any positive continuously differentiable function and $t_{2}$ is a smooth function of $b^{2}$.

\end{thm}
\textbf{Remark:} 1) we assume that $\b$ is closed and conformal with respect to $\a$, i.e. (\ref{bcij}) holds. According to the relate discussions for isotropic Berwald metrics  \cite{Chen-Shen, Guo-Liu-Mo, Tayebi-Rafie, Tayebi-Najafi},  we believe that the assumption here is reasonable and appropriate.

2) It should be pointed out that if the scalar function $c(x)=0$, then according to Proposition \ref{Bewald curvature}, $B_{j}{}^{i}{}_{kl}=0$, namely, $F=\gab$ is a Berwald metric for any function $\phi(b^{2},s)$. So it will be regarded as a trivial case.

\section{Preliminaries}

Let $F$ be a Finsler metric on an $n$-dimensional manifold $M$ and $G^{i}$ be the geodesic coefficients of $F$, which are defined by
\begin{eqnarray*}
G^{i}=\frac{1}{4}g^{il}\left\{[F^{2}]_{x^{k}y^{l}}y^{k}-[F^{2}]_{x^{l}}\right\},
\end{eqnarray*}
where $(g^{ij}):=\left(\frac{1}{2}[F^{2}]_{y^{i}y^{j}}\right)^{-1}$. For a Riemannian metric, the spray coefficients are determined by its Christoffel symbols as $G^{i}(x,y)=\frac{1}{2}\Gamma^{i}_{jk}(x)y^{j}y^{k}$.

\begin{lem}\cite{yct-zhm-onan}\label{Finsler}
 Let $F=\gab$ be a  general $\ab$-metric on an $n$-dimensional manifold $M$.
Then the function $F$ is a regular Finsler metric for any Riemannian metric $\alpha$ and  any $1$-form $\beta$ if and only if $\phi(b^{2},s)$ is a positive smooth function defined on the domain $|s|\leq b<b_o$ for some positive number (maybe infinity) $b_o$ satisfying
\begin{eqnarray}\label{ppp}
\p-s\pt>0,\quad\p-s\pt+(b^2-s^2)\ptt>0,
\end{eqnarray}
when $n\geq 3$ or
\begin{eqnarray}\label{ppp1}
\p-s\pt+(b^2-s^2)\ptt>0,
\end{eqnarray}
when $n=2$.
\end{lem}

Let $\alpha=\sqrt{a_{ij}(x)y^iy^j}$  and $\beta= b_i(x)y^i$.
Denote the coefficients of the covariant derivative of
$\b$ with respect to $\a$ by $b_{i|j}$, and let
\begin{eqnarray*}
&r_{ij}=\frac{1}{2}(b_{i|j}+b_{j|i}),~s_{ij}=\frac{1}{2}(b_{i|j}-b_{j|i}),
~r_{00}=r_{ij}y^iy^j,~s^i{}_0=a^{ij}s_{jk}y^k,&\\
&r_i=b^jr_{ji},~s_i=b^js_{ji},~r_0=r_iy^i,~s_0=s_iy^i,~r^i=a^{ij}r_j,~s^i=a^{ij}s_j,~r=b^ir_i,&
\end{eqnarray*}
where $(a^{ij}):=(a_{ij})^{-1}$ and $b^{i}:=a^{ij}b_{j}$. It is easy to see that $\b$ is closed if and only if $s_{ij}=0$.

\begin{lem}\cite{yct-zhm-onan}
 the spray coefficients $G^i$ of a general $(\alpha,\beta)$-metric $F=\gab$ are related to the spray coefficients ${}^\a G^i$ of
$\a$ and given by
\begin{eqnarray}\label{Gi}
G^i&=&{}^\a G^i+\a Q s^i{}_0+\left\{\Theta(-2\a Q s_0+r_{00}+2\a^2
R r)+\a\Omega(r_0+s_0)\right\}\frac{y^i}{\a}\nonumber\\
&&+\left\{\Psi(-2\a Q s_0+r_{00}+2\a^2 R
r)+\a\Pi(r_0+s_0)\right\}b^i -\a^2 R(r^i+s^i),
\end{eqnarray}
where
$$Q=\frac{\pt}{\p-s\pt},\quad R=\frac{\po}{\p-s\pt},$$
$$\Theta=\frac{(\p-s\pt)\pt-s\p\ptt}{2\p\big(\p-s\pt+(b^2-s^2)\ptt\big)},
\quad\Psi=\frac{\ptt}{2\big(\p-s\pt+(b^2-s^2)\ptt\big)},$$
$$\Pi=\frac{(\p-s\pt)\pot-s\po\ptt}{(\p-s\pt)\big(\p-s\pt+(b^2-s^2)\ptt\big)},\quad
\Omega=\frac{2\po}{\p}-\frac{s\p+(b^2-s^2)\pt}{\p}\Pi.$$
\end{lem}
Note that $\phi_1$ means the derivation of $\phi$ with respect to the first variable $b^2$. In the following, we will introduce an important non-Riemannian quantity.
\begin{definition} \cite{szm-diff}
 Let
\begin{eqnarray}\label{DT}
 B_{j}{}^{i}{}_{kl}: = \frac{\partial^{3}G^{i}}{\partial y^{j}\partial y^{k} \partial y^{l}},
\end{eqnarray}
where $G^{i}$ are the spray coefficients of $F$. The tensor
$B:=B_{j}{}^{i}{}_{kl} \partial_{i}\otimes dx^{j}\otimes dx^{k}\otimes
dx^{l}$ is called {\em Berwald  tensor}. A Finsler metric is called
{\it a Berwald metric} if the Berwald tensor vanishes, i.e. the spray coefficients $G^{i}=G^{i}(x,y)$ are quadratic in $y\in T_{x}M$ at every point $x\in M$.
\end{definition}

\section{Berwald curvature of general $(\alpha,\beta)$-metrics}
In this section, we will compute the Berwald curvature of a general $\ab$-metric.
\begin{prop}\label{Bewald curvature}
Let $F=\gab$ be a general $\ab$-metric on an $n$-dimensional manifold $M$. Suppose that $\b$ satisfies (\ref{bcij}), then the Berwald curvature of $F$ is given by
\begin{eqnarray}
B_{j}{}^{i}{}_{kl}&=&\frac{c}{\alpha}\left\{\left[(E-sE_{2})a_{kl}+E_{22}b_{l}b_{k}\right]\delta^{i}{}_{j}+\frac{1}{\alpha^{2}}\left[\frac{s}{\alpha}(3E_{22}+sE_{222})y_{l}y_{j}-
(E_{22}+sE_{222})b_{l}y_{j}\right]b_{k}y^{i}\right\}(k\rightarrow l\rightarrow j\rightarrow k)\nonumber\\
&&-\frac{c}{\alpha^{2}}\left\{sE_{22}\left[(y_{k}b_{l}+y_{l}b_{k})\delta^{i}{}_{j}+a_{jl}b_{k}y^{i}\right]+\frac{1}{\alpha}(E-sE_{2}-s^{2}E_{22})
(y_{l}\delta^{i}{}_{j}+a_{lj}y^{i})y_{k}\right\}(k\rightarrow l\rightarrow j\rightarrow k)\nonumber\\
&&+\frac{c}{\alpha^{2}}\left[\frac{1}{\alpha^{3}}(3E-3sE_{2}-6s^{2}E_{22}-s^{3}E_{222})y_{k}y_{j}y_{l}+E_{222}b_{l}b_{k}b_{j}\right]y^{i}\nonumber\\
&&+\frac{c}{\alpha}\left[(H_{2}-sH_{22})(b_{j}-\frac{s}{\alpha}y_{j})a_{kl}-\frac{1}{\alpha^{2}}(H_{2}-sH_{22}-s^{2}H_{222})b_{l}y_{j}y_{k}
-\frac{sH_{222}}{\alpha}b_{k}b_{l}y_{j}\right]b^{i}(k\rightarrow l\rightarrow j\rightarrow k)\nonumber\\
&&+\frac{c}{\alpha}\left[\frac{s}{\alpha^{3}}(3H_{2}-3sH_{22}-s^{2}H_{222})y_{j}y_{k}y_{l}+H_{222}b_{l}b_{k}b_{j}\right]b^{i},\label{BC}
\end{eqnarray}
where $y_{i}:=a_{ij}y^{j}$ and $b^{i}:=a^{ij}b_{j}$, $c=c(x)\not=0$ is a scalar function on $M$.
\begin{eqnarray}
 E:&=&\frac{\pt+2s\po}{2\p}
-H\frac{s\p+(b^2-s^2)\pt}
{\p},\label{T}\\
H:&=&\frac{\ptt-2(\po-s\pot)}{2\big[\p-s\pt+(b^2-s^2)\ptt\big]}.\label{H}
\end{eqnarray}
\end{prop}
\begin{proof}
By (\ref{bcij}), we have
\begin{eqnarray}\label{tu}
r_{00}=c\a^2,r_0=c\b,r=cb^2,r^i=cb^i,s^i{}_0=0,s_0=0,s^i=0.
\end{eqnarray}
Substituting (\ref{tu}) into (\ref{Gi}) yields
\begin{eqnarray}
G^i&=&{}^\a G^i+c\a\left\{\Theta(1+2Rb^2)+s\Omega\right\}y^i+c\a^2\left\{\Psi(1+2Rb^2)+s\Pi-R\right\}b^i\nonumber\\
&=&{}^\a G^i+c\a E y^i+c\a^2H b^i,\label{gp}
\end{eqnarray}
where $E$ and $H$ are given by (\ref{T}) and (\ref{H}) respectively.
Note that
\begin{eqnarray}\label{aby}
\alpha_{y^{i}}=\frac{y_{i}}{\alpha},~~ s_{y^{i}}=\frac{\alpha b_{i}-sy_{i}}{\a^{2}},
\end{eqnarray}
where $y_{i}:=a_{ij}y^{j}$.
Put
\begin{eqnarray}\label{Wi}
W^{i}:=\a Ey^{i}+\alpha^{2} H b^{i}.
\end{eqnarray}
Differentiating (\ref{Wi}) with respect to $y^{j}$ yields
\begin{eqnarray}\label{Wij}
\frac{\partial W^{i}}{\partial y^{j}}=\a E\delta^{i}{}_{j}+(E\a_{y^{j}}+\alpha E_{2}s_{y^{j}})y^{i}+\big\{[\alpha^{2}]_{y^{j}} E+\a^{2}E_{2}s_{y^{j}}\big\} b^{i}.
\end{eqnarray}
Differentiating (\ref{Wij}) with respect to $y^{k}$ yields
\begin{eqnarray}\label{Wijk11}
\frac{\partial^{2} W^{i}}{\partial y^{j}\partial y^{k}}&=&\big[(E\a_{y^{k}}+\a E_{2}s_{y^{k}})\delta^{i}{}_{j}+E_{2}s_{y^{k}}\a_{y^{j}}y^{i}+H_{2}[\a^{2}]_{y^{j}}s_{y^{k}}b^{i}\big](k\leftrightarrow j)\nonumber\\
&&+\big(E\a_{y^{j}y^{k}}+\a E_{22}s_{y^{k}}s_{y^{j}}+\a E_{2}s_{y^{j}y^{k}}\big)y^{i}\nonumber\\
&&+\big\{[\a^{2}]_{y^{j}y^{k}}H+\a^{2}H_{22}s_{y^{k}}s_{y^{j}}+\a^{2}H_{2}s_{y^{j}y^{k}}\big\}b^{i},
\end{eqnarray}
where $k\leftrightarrow j$ denotes symmetrization. Therefore, it follows from (\ref{Wijk11}) that
\begin{eqnarray}\label{Wijk}
\frac{\partial^{3} W^{i}}{\partial y^{j}\partial y^{k}\partial y^{l}}&=& \big[E_{2}(\a_{y^{k}}s_{y^{l}}+\a_{y^{l}}s_{y^{k}}+\a s_{y^{k}y^{l}})+E\a_{y^{k}y^{l}}+\a E_{22}s_{y^{l}}s_{y^{k}}\big]\delta^{i}{}_{j}(k\rightarrow l\rightarrow j\rightarrow k)\nonumber\\
&&+\big[E_{2}(s_{y^{k}}\a_{y^{j}y^{l}}+\a_{y^{k}}s_{y^{j}y^{l}})+E_{22}(\a_{y^{k}}s_{y^{j}}+\a s_{y^{k}y^{j}})s_{y^{l}}\big]y^{i}
(k\rightarrow l\rightarrow j\rightarrow k)\nonumber\\
&&+\Large\{H_{2}\big([\a^{2}]_{y^{k}y^{l}}s_{y^{j}}+[\a^{2}]_{y^{k}}s_{y^{j}y^{l}}\big)+H_{22}\big([\a^{2}]_{y^{k}}s_{y^{l}}s_{y^{j}}
+\a^{2}s_{y^{k}y^{l}}s_{y^{j}}\big)\Large\}b^{i}(k\rightarrow l\rightarrow j\rightarrow k)\nonumber\\
&&+\big(E\a_{y^{j}y^{k}y^{l}}+\a E_{222}s_{y^{j}}s_{y^{k}}s_{y^{l}}+\a E_{2}s_{y^{j}y^{k}y^{l}}\big)y^{i}\nonumber\\
&&+\big\{H[\a^{2}]_{y^{j}y^{k}y^{l}}+\a^{2} H_{222}s_{y^{j}}s_{y^{k}}s_{y^{l}}+\a^{2} H_{2}s_{y^{j}y^{k}y^{l}}\big\}b^{i},
\end{eqnarray}
where $k\rightarrow l\rightarrow j\rightarrow k$ denotes cyclic permutation.
It follows from (\ref{aby}) that
\begin{eqnarray}\label{asd}
[\a^{2}]_{y^{l}}&=& 2y_{l},~~~ [\a^{2}]_{y^{l}y^{j}}=2a_{lj},~~~ [\a^{2}]_{y^{l}y^{j}y^{k}}=0,\label{a2d} \\
\a_{y^{l}y^{j}}&=& \frac{1}{\a}\big(a_{lj}-\frac{y_{l}}{\a}\frac{y_{j}}{\a}\big),~~ \a_{y^{l}y^{j}y^{k}}=-\frac{1}{\a^{3}}[a_{kl}y_{j}(k\rightarrow l\rightarrow j\rightarrow k)-\frac{3}{\a^{2}}y_{l}y_{j}y_{k}],\label{ad}\\
s_{y^{l}y^{j}}&=& -\frac{1}{\a^{2}}\big[s a_{lj}+\frac{1}{\a}(b_{l}y_{j}+b_{j}y_{l})-\frac{3s}{\a^{2}}y_{l}y_{j}\big],\label{s2d}\\
s_{y^{l}y^{j}y^{k}}&=& \frac{1}{\a^{5}}\big\{[\a(3s y_{j}-\a b_{j})a_{lk}+3b_{k}y_{l}y_{j}](k\rightarrow l\rightarrow j\rightarrow k)-\frac{15s}{\a}y_{k}y_{l}y_{j}\big\}.\label{s3d}
\end{eqnarray}
Plugging (\ref{a2d})-(\ref{s3d}) into (\ref{Wijk}) yields
\begin{eqnarray}
\frac{\partial ^{3}W^{i}}{\partial y^{j}\partial y^{k}\partial y^{l}}&=&\frac{1}{\a}\left\{\big[(E-sE_{2})a_{kl}+E_{22}b_{l}b_{k}\big]\delta^{i}{}_{j}+
\frac{1}{\a^{2}}\big[\frac{s}{\a}(3E_{22}+sE_{222})y_{l}-(E_{22}+sE_{222})b_{l}\big]y_{j}b_{k}y^{i}\right\}(k\rightarrow l\rightarrow j\rightarrow k)\nonumber\\
&&-\frac{1}{\a^{2}}\left\{sE_{22}\big[(y_{k}b_{l}+y_{l}b_{k})\delta^{i}{}_{j}+a_{jl}b_{k}y^{i}\big]
+\frac{1}{\a}(E-sE_{2}-s^{2}E_{22})(y_{l}\delta^{i}_{j}+a_{jl}y^{i})y_{k}\right\}(k\rightarrow l\rightarrow j\rightarrow k)\nonumber\\
&&+\frac{1}{\alpha^{2}}\left[\frac{1}{\alpha^{3}}(3E-3sE_{2}-6s^{2}E_{22}-s^{3}E_{222})y_{k}y_{j}y_{l}+E_{222}b_{l}b_{k}b_{j}\right]y^{i}\nonumber\\
&&+\frac{1}{\alpha}\left[(H_{2}-sH_{22})(b_{j}-\frac{s}{\alpha}y_{j})a_{kl}-\frac{1}{\alpha^{2}}(H_{2}-sH_{22}-s^{2}H_{222})b_{l}y_{j}y_{k}
-\frac{sH_{222}}{\alpha}b_{k}b_{l}y_{j}\right]b^{i}(k\rightarrow l\rightarrow j\rightarrow k)\nonumber\\
&&+\frac{1}{\alpha}\left[\frac{s}{\alpha^{3}}(3H_{2}-3sH_{22}-s^{2}H_{222})y_{j}y_{k}y_{l}+H_{222}b_{l}b_{k}b_{j}\right]b^{i}.\label{Wijk1}
\end{eqnarray}
It follows from ${}^\a G^i(x,y)=\frac{1}{2}\Gamma^{i}_{jk}(x)y^{j}y^{k}$ that
\begin{eqnarray}\label{ajlk}
\frac{\partial ^{3}{}^\a G^{i}}{\partial y^{j}\partial y^{k}\partial y^{l}}=0.
\end{eqnarray}
By (\ref{DT}), (\ref{gp}), (\ref{Wi}), (\ref{Wijk1}) and (\ref{ajlk}), we obtain (\ref{BC}).

\end{proof}

\section{General $\ab$-metrics with isotropic Berwald curvature}
In this section, we will classify general $\ab$-metrics with isotropic Berwald curvature under certain condition. Firstly, we show the following

\begin{lem}\label{HE}
Suppose that $\b$ satisfies (\ref{bcij}), then a general $\ab$-metric $F=\gab$ is of isotropic Berwald curvature if and only if
\begin{eqnarray}
E-sE_{2}-\rho(x)(\phi-s\phi_{2})=0, \label{lresult1}\\
H_{2}-sH_{22}=0, \label{lresult2}
 \end{eqnarray}
where $\rho(x)=\frac{\tau(x)}{c(x)}$, $E$ and $H$ are given by (\ref{T}) and (\ref{H}) respectively. In particular, $F$ is a Berwald metric if and only if $E-sE_{2}=H_{2}-sH_{22}=0$ .
\end{lem}
\begin{proof}
For a general $\ab$-metric $F=\a\phi(b^{2},s)$, a direct computation yields
\begin{eqnarray}
F_{y^{j}}&=& \a_{y^{j}}\phi+\a\phi_{2}s_{y^{j}},\nonumber\\
F_{y^{j}y^{k}} &=& \a_{y^{j}y^{k}}\phi+(\alpha_{y^{j}} s_{y^{k}}+\alpha_{y^{k}}s_{y^{j}})\phi_{2}+\alpha\phi_{22}s_{y^{k}}s_{y^{j}}+\alpha\phi_{2}s_{y^{j}y^{k}},\label{F2}\\
F_{y^{j}y^{k}y^{l}} &=& [(\a_{y^{j}y^{k}}s_{y^{l}}+\a_{y^{j}}s_{y^{k}y^{l}})\phi_{2}+(\a_{y^{j}}s_{y^{l}}+\a s_{y^{j}y^{l}})s_{y^{k}}\phi_{22}](j\rightarrow k\rightarrow l\rightarrow j)\nonumber\\
&&+\a_{y^{j}y^{k}y^{l}}\phi+\alpha \phi_{222}s_{y^{l}}s_{y^{k}}s_{y^{j}}+\alpha\phi_{2}s_{y^{j}y^{k}y^{l}},\label{F3}
\end{eqnarray}
Plugging (\ref{aby}) and (\ref{ad})-(\ref{s3d}) into (\ref{F2}) and (\ref{F3}) yields
\begin{eqnarray}
F_{y^{j}y^{k}} &=& \frac{1}{\a}(\phi-s\phi_{2})a_{jk}-\frac{s\phi_{22}}{\a^{2}}(b_{k}y_{j}+b_{j}y_{k})+\frac{\phi_{22}}{\alpha}b_{j}b_{k}
-\frac{1}{\alpha^{3}}(\phi-s\phi_{2}-s^{2}\phi_{22})y_{j}y_{k},\label{F22}\\
F_{y^{j}y^{k}y^{l}} &=& \frac{1}{\a^{2}}\big[\frac{1}{\a}(s\phi_{2}+s^{2}\phi_{22}-\phi)a_{kl}y_{j}+\frac{s}{\alpha^{2}}(3\phi_{22}+s\phi_{222})b_{k}y_{l}y_{j}
-s\phi_{22}a_{jl}b_{k}\nonumber\\
&&-\frac{1}{\a}(\phi_{22}+s\phi_{222})b_{l}b_{j}y_{k}\big](j\rightarrow k\rightarrow l\rightarrow j)\nonumber\\
&&+\frac{1}{\a^{5}}(\phi-3s\phi_{2}-6s^{2}\phi_{22}-s^{3}\phi_{222})y_{j}y_{k}y_{l}+\frac{1}{\a^{2}}\phi_{222}b_{l}b_{k}b_{j}.\label{F222}
\end{eqnarray}
Suppose $F$ be of isotropic Berwald curvature, by (\ref{isotropic}), (\ref{F22}) and (\ref{F222}), we obtain
\begin{eqnarray}
B_{j}{}^{i}{}_{kl}&=&\frac{\tau(x)}{\a}\big[(\phi-s\phi_{2})a_{jk}-\frac{s\phi_{22}}{\a}(b_{k}y_{j}+b_{j}y_{k})+\phi_{22}b_{j}b_{k}
-\frac{1}{\alpha^{2}}(\phi-s\phi_{2}-s^{2}\phi_{22})y_{j}y_{k}\big]\delta^{i}{}_{l}(j\rightarrow k\rightarrow l\rightarrow j)\nonumber\\
&&+\frac{\tau(x)}{\a^{2}}\big[\frac{1}{\a}(s\phi_{2}+s^{2}\phi_{22}-\phi)a_{kl}y_{j}+\frac{s}{\alpha^{2}}(3\phi_{22}+s\phi_{222})b_{k}y_{l}y_{j}
-s\phi_{22}a_{jl}b_{k}\nonumber\\
&&-\frac{1}{\a}(\phi_{22}+s\phi_{222})b_{l}b_{j}y_{k}\big]y^{i}(j\rightarrow k\rightarrow l\rightarrow j)\nonumber\\
&&+\frac{\tau(x)}{\a^{2}}\big[\frac{1}{\a^{3}}(3\phi-3s\phi_{2}-6s^{2}\phi_{22}-s^{3}\phi_{222})y_{j}y_{k}y_{l}
+\phi_{222}b_{l}b_{k}b_{j}\big]y^{i}.\label{isotropic1}
\end{eqnarray}
By (\ref{BC}) and (\ref{isotropic1}), we obtain
\begin{eqnarray}\label{spara}
T_{1}+\a T_{2}=0,
\end{eqnarray}
where
\begin{eqnarray}
T_{1}:&=& \a^{4}\big\{[E-sE_{2}-\rho(\phi-s\phi_{2})]a_{kl}\delta^{i}_{j}+(E_{22}-\rho\phi_{22})b_{l}b_{k}\delta^{i}_{j}+(H_{2}-sH_{22})a_{kl}b_{j}\big\}(j\rightarrow k\rightarrow l\rightarrow j)\nonumber\\
&&-\a^{2}
[E-sE_{2}-s^{2}E_{22}-\rho(\phi-s\phi_{2}-s^{2}\phi_{22})](y_{l}\delta^{i}_{j}+a_{lj}y^{i})y_{k}(j\rightarrow k\rightarrow l\rightarrow j)+\a^{4}H_{222}b_{l}b_{k}b_{j}b^{i}\nonumber\\
&&-\a^{2}\big\{[E_{22}+sE_{222}-\rho(\phi_{22}+s\phi_{222})]b_{k}y^{i}+(H_{2}-sH_{22}-s^{2}H_{222})y_{k}b^{i})
\big\}b_{l}y_{j}(j\rightarrow k\rightarrow l\rightarrow j)\nonumber\\
&&+[3E-3sE_{2}-6s^{2}E_{22}-s^{3}E_{222}-\rho(3\phi-3s\phi_{2}-6s^{2}\phi_{22}-s^{3}\phi_{222})]y_{j}y_{k}y_{l}y^{i},\label{T1}\\
T_{2}:&=& -\a^{2} s\big\{(E_{22}-\rho\phi_{22})[(y_{k}b_{l}+y_{l}b_{k})\delta^{i}{}_{j}+a_{jl}b_{k}y^{i}]+[(H_{2}-sH_{22})y_{j}a_{kl}+H_{222}b_{k}b_{l}y_{j}]b^{i}\big\}
(j\rightarrow k\rightarrow l\rightarrow j)\nonumber\\
&&+s\big[3E_{22}+sE_{222}-\rho(3\phi_{22}+s\phi_{222})\big]y_{l}y_{j}b_{k}y^{i}(j\rightarrow k\rightarrow l\rightarrow j)\nonumber\\
&&+\a^{2}(E_{222}-\rho\phi_{222})b_{l}b_{k}b_{j}y^{i}+s(3H_{2}-3sH_{22}-s^{2}H_{222})y_{j}y_{k}y_{l}b^{i}.\label{T2}
\end{eqnarray}
By (\ref{spara}), we know that
$$
T_{1}=0,~~~ T_{2}=0.
$$
For $s\neq 0$, it follows from $T_{2}y^{j}y^{k}y^{l}=0$ that
\begin{eqnarray}\label{yjykyl}
(E_{222}-\rho\phi_{222})y^{i}-\a H_{222}b^{i}=0.
\end{eqnarray}
 Both rational part and irrational part of (\ref{yjykyl}) equal zero, namely
\begin{eqnarray}\label{eq1}
E_{222}-\rho\phi_{222}=0, ~~~ H_{222}=0.
\end{eqnarray}
Plugging (\ref{eq1}) into $T_{2}=0$ yields
\begin{eqnarray}\label{eq2}
&&-\a^{2} \big\{(E_{22}-\rho\phi_{22})[(y_{k}b_{l}+y_{l}b_{k})\delta^{i}{}_{j}+a_{jl}b_{k}y^{i}]+(H_{2}-sH_{22})a_{kl}y_{j}b^{i}\big\}
(j\rightarrow k\rightarrow l\rightarrow j)\nonumber\\
&&+3(E_{22}-\rho\phi_{22})y_{l}y_{j}b_{k}y^{i}(j\rightarrow k\rightarrow l\rightarrow j)+3(H_{2}-sH_{22})y_{j}y_{k}y_{l}b^{i}=0.
\end{eqnarray}
Contracting (\ref{eq2}) by $b^{j}b^{k}b^{l}$ yields
\begin{eqnarray}\label{eq3}
(E_{22}-\rho\phi_{22})b^{2}(b^{2}-3s^{2})y^{i}+\a s[2b^{2}(E_{22}-\rho\phi_{22})+(H_{2}-sH_{22})(b^{2}-s^{2})]b^{i}=0.
\end{eqnarray}
Hence, it follows from (\ref{eq3}) that
\begin{eqnarray}\label{eq4}
E_{22}-\rho\phi_{22}=0, ~~~~ H_{2}-sH_{22}=0.
\end{eqnarray}
Inserting (\ref{eq1}) and (\ref{eq4}) into $T_{1}=0$ yields
\begin{eqnarray}\label{eq5}
&&[E-sE_{2}-\rho(\phi-s\phi_{2})]\big\{\a^{2}[\a^{2}a_{kl}\delta^{i}_{j}-(y_{l}\delta^{i}_{j}+a_{lj}y^{i})y_{k}](j\rightarrow k\rightarrow l\rightarrow j)+3y_{j}y_{k}y_{l}y^{i}\big\}=0.
\end{eqnarray}
Multiplying (\ref{eq5}) by $b^{j}b^{k}b^{l}$ yields
\begin{eqnarray}\label{eq6}
[E-sE_{2}-\rho(\phi-s\phi_{2})](b^{2}-s^{2})(\a b^{i}-s y^{i})=0.
\end{eqnarray}
Hence, it is easy to see from (\ref{eq6}) that
\begin{eqnarray}\label{eq7}
E-sE_{2}-\rho(\phi-s\phi_{2})=0.
\end{eqnarray}
Note that
\begin{eqnarray*}
E_{222}-\rho\phi_{222}&=& (E_{22}-\rho\phi_{22})_{2},~~
s(E_{22}-\rho\phi_{22})=-\big[E-sE_{2}-\rho(\phi-s\phi_{2})\big]_{2},~~ sH_{222}=-(H_{2}-sH_{22})_{2}.
\end{eqnarray*}
Therefore, (\ref{eq7}) implies that the first equalities of (\ref{eq1}) and (\ref{eq4}) hold. The second equality of (\ref{eq4}) implies that the second equality of (\ref{eq1}) holds. Moreover,  if a general $\ab$-metric $F=\a\phi(b^{2},s)$ is of isotropic Berwald curvature, then
(\ref{lresult1}) and (\ref{lresult2}) hold.

Conversely, suppose that (\ref{lresult1}) and (\ref{lresult2}) hold. Setting $\psi:=E-\rho(x)\phi$, then (\ref{lresult1}) is equivalent to
\begin{eqnarray}\label{eq8}
\psi-s\psi_{2}=0.
\end{eqnarray}
By solving Eq. (\ref{eq8}), we obtain $\psi=\frac{1}{2}\sigma(b^{2})s$. Hence,
\begin{eqnarray}\label{eq9}
E-\rho(x)\phi=\frac{1}{2}\sigma(b^{2})s.
\end{eqnarray}
By (\ref{lresult2}), there exist two functions $t_{1}(b^{2})$ and $t_{2}(b^{2})$ such that
\begin{eqnarray}\label{eq10}
H=\frac{1}{2}\big[t_{1}(b^{2})+t_{2}(b^{2})s^{2}\big].
\end{eqnarray}
By (\ref{gp}), (\ref{eq9}) and (\ref{eq10}),
\begin{eqnarray}\label{eq11}
G^{i}-\tau(x)Fy^{i}&=&{}^\a G^i+c(x)\a E y^i+c(x)\a^2H b^i-\tau(x)Fy^{i}\nonumber\\
&=&{}^\a G^i+\a\big[c(x)E-\tau(x)\phi\big]y^{i}+\frac{1}{2}c(x)\a^{2}\big[t_{1}(b^{2})+t_{2}(b^{2})s^{2}\big]b^{i}\nonumber\\
&=& {}^\a G^i+\frac{1}{2}c(x)\big\{\sigma(b^{2})\b y^{i}+[t_{1}(b^{2})\a^{2}+t_{2}(b^{2})\b^{2}]b^{i}\big\}.
\end{eqnarray}
Hence, $G^{i}$ are quadratic in $y=y^{i}\frac{\partial}{\partial x^{i}}|_{x}$. On the other hand, it follows from (\ref{DT}) that
\begin{eqnarray}\label{eq12}
(G^{i}-\tau(x)Fy^{i})_{y^{j}y^{k}y^{l}}=B_{j}{}^{i}{}_{kl}-\tau(x)(F_{y^{j}y^{k}}\delta^{i}{}_{l}+F_{y^{j}y^{l}}\delta^{i}{}_{k}
+F_{y^{l}y^{k}}\delta^{i}{}_{j}+F_{y^{j}y^{k}y^{l}}y^{i}).
\end{eqnarray}
Using (\ref{eq11}) and (\ref{eq12}), we obtain that $F=\a\phi(b^{2},s)$ is of isotropic Berwald curvature.

Observe that $F$ is a Berwald metric if and only if its Berwald curvature $B_{j}{}^{i}{}_{kl}=0$. By (\ref{isotropic}), we obtain that $F$ is a Berwald metric if and only if $\tau(x)=0$, i.e., $\rho(x)=0$. Hence, by the above process of proof, we get that $F$ is a Berwald metric if and only if $E-sE_{2}=H_{2}-sH_{22}=0$ .
\end{proof}

\begin{proof}[Proof of Theorem \ref{main2}]
Suppose that $\b$ satisfies (\ref{bcij}), By Lemma \ref{HE},  a general $\ab$-metric $F=\gab$ is of isotropic Berwald curvature if and only if
(\ref{lresult1}) and (\ref{lresult2}) hold. Then there exist three functions $\sigma(b^{2})$, $t_{1}(b^{2})$ and $t_{2}(b^{2})$ such that (\ref{eq9})
and (\ref{eq10}) hold. Plugging (\ref{eq9})-(\ref{eq10}) into (\ref{T}) and (\ref{H}) yields
\begin{eqnarray}
\frac{\pt+2s\po}{2\p}
-\frac{1}{2}(t_{1}+t_{2}s^{2})\frac{s\p+(b^2-s^2)\pt}
{\p}&=&\rho\phi+\frac{1}{2}\sigma s,\label{TTT}\\
\frac{\ptt-2(\po-s\pot)}{2\big[\p-s\pt+(b^2-s^2)\ptt\big]}&=& \frac{1}{2}\big(t_{1}+t_{2}s^{2}\big),\label{HHH}
\end{eqnarray}
where we use $\rho$, $\sigma$, $t_{1}$ and $t_{2}$ instead of $\rho(x)$, $\sigma(b^{2})$, $t_{1}(b^{2})$ and $t_{2}(b^{2})$, respectively.
(\ref{TTT}) and (\ref{HHH}) are equivalent to
\begin{eqnarray}
\big[1-(t_{1}+t_{2}s^{2})(b^2-s^2)\big]\phi_{2}+2s\phi_{1}-s\big[(t_{1}+t_{2}s^{2})+\sigma\big]\phi-2\rho\phi^{2}&=&0, \label{TTT1}\\
\big[1-(b^{2}-s^{2})(t_{1}+t_{2}s^{2})\big]\phi_{22}-2\phi_{1}+2s\phi_{12}+s(t_{1}+t_{2}s^{2})\phi_{2}-(t_{1}+t_{2}s^{2})\phi &=& 0. \label{HHH1}
\end{eqnarray}
Differentiating (\ref{TTT1}) with respect to $s$ yields
\begin{eqnarray}\label{Ds}
&&\big[1-(t_{1}+t_{2}s^{2})(b^2-s^2)\big]\phi_{22}+2\phi_{1}+2s\phi_{12}+s\big(t_{1}-\sigma-2b^{2}t_{2}+3t_{2}s^{2}\big)\phi_{2}\nonumber\\
&&-\big(t_{1}+\sigma+3t_{2}s^{2}\big)\phi-4\rho\phi\phi_{2}=0.
\end{eqnarray}
From $(\ref{Ds})-(\ref{HHH1})$, we obtain
\begin{eqnarray}\label{Ds1}
4\phi_{1}-s\big[2t_{2}(b^{2}-s^{2})+\sigma\big]\phi_{2}-(\sigma+2t_{2}s^{2})\phi-4\rho\phi\phi_{2}=0.
\end{eqnarray}
From $(\ref{TTT1})\times2-(\ref{Ds1})\times s$, we obtain
\begin{eqnarray}\label{Ds2}
\big[2-2t_{1}(b^{2}-s^{2})+\sigma s^{2}\big]\phi_{2}-(2t_{1}+\sigma)s\phi-4\rho\phi^{2}+4\rho s\phi\phi_{2}=0.
\end{eqnarray}
Note that (\ref{Ds2}) is equivalent to
\begin{eqnarray}\label{Ds3}
\big(\frac{2-2t_{1}(b^{2}-s^{2})+\sigma s^{2}}{\phi^{2}}\big)_{2}+\big(\frac{8\rho s}{\phi}\big)_{2}=0.
\end{eqnarray}
\textbf{Case 1 $\rho\neq 0$}

\textbf{1) $2-2t_{1}(b^{2}-s^{2})+\sigma s^{2}\neq0$}

Integrating (\ref{Ds3}) yields
\begin{eqnarray}
\phi=\frac{4\rho s+\sqrt{2k(1-t_{1}b^{2})+(16\rho^{2}+2kt_{1}+k\sigma)s^{2}}}{k},
\end{eqnarray}
where $k=k(b^{2})$ is any non-zero smooth function. Then the corresponding general $\ab$-metric $F=\a\phi(b^{2},s)$ is a Randers metric.

\textbf{2) $2-2t_{1}(b^{2}-s^{2})+\sigma s^{2}=0$}

In this case, (\ref{Ds3}) is reduced to $(\frac{s}{\phi})_{2}=0$, it is easy to obtain that $\phi=\frac{s}{a(b^{2})}$. Hence, the corresponding general $\ab$-metric $F=\a\phi(b^{2},s)$ is a Kropina metrics.\\
\textbf{Case 2 $\rho = 0$}

In this case, it follows from (\ref{isotropic}) that $F$ is a Berwald metric. (\ref{Ds2}) is reduced to
\begin{eqnarray}\label{Dss2}
\big[2-2t_{1}(b^{2}-s^{2})+\sigma s^{2}\big]\phi_{2}-(2t_{1}+\sigma)s\phi=0.
\end{eqnarray}

\textbf{i) $2-2t_{1}(b^{2}-s^{2})+\sigma s^{2}\neq0$}

By (\ref{Dss2}), we obtain
\begin{eqnarray}\label{solutions11}
\phi= t_{3}(b^{2})\sqrt{2(1-b^{2}t_{1})+(\sigma+2t_{1})s^{2}},
\end{eqnarray}
where $t_{3}(b^{2})$ is any positive smooth function. Hence, in this case, the corresponding  general $\ab$-metric
$F=\a\phi(b^{2}, s)$
is a Riemannian metric.

\textbf{ii) $2-2t_{1}(b^{2}-s^{2})+\sigma s^{2}=0$}

Note that $\phi>0$ and $s\neq0$. In this case, (\ref{Dss2}) is equivalent to
\begin{eqnarray}\label{TT6}
\sigma+2t_{1}=0, ~~ 2-2(b^{2}-s^{2})t_{1}+\sigma s^{2}=0.
\end{eqnarray}
By (\ref{TT6}), we have
\begin{eqnarray}\label{f1f2}
\sigma=-\frac{2}{b^{2}}, ~~ t_{1}=\frac{1}{b^{2}}.
\end{eqnarray}
In this case, (\ref{Ds1}) imply that (\ref{TTT1}) holds. By the above caculations, it is easy to see that (\ref{TTT1}) and (\ref{Ds1}) imply (\ref{HHH1}). Therefore, we only need to solve (\ref{Ds1}).  Plugging (\ref{f1f2}) into (\ref{Ds1}) yields
\begin{eqnarray}\label{1order}
\phi_{1}+\frac{1}{2}s[\frac{1}{b^{2}}-(b^{2}-s^{2})t_{2}]\phi_{2}=\frac{1}{2}(-\frac{1}{b^{2}}+t_{2}s^{2})\phi.
\end{eqnarray}
The characteristic equation of PDE (\ref{1order}) is
\begin{eqnarray}\label{chara}
\frac{db^{2}}{1}=\frac{ds}{\frac{1}{2}s[\frac{1}{b^{2}}-(b^{2}-s^{2})t_{2}]}=\frac{d\phi}{\frac{1}{2}(-\frac{1}{b^{2}}+t_{2}s^{2})\phi}.
\end{eqnarray}
Firstly, we solve
\begin{eqnarray}\label{chara1}
\frac{db^{2}}{1}=\frac{ds}{\frac{1}{2}s[\frac{1}{b^{2}}-(b^{2}-s^{2})t_{2}]}.
\end{eqnarray}
(\ref{chara1}) is equivalent to
\begin{eqnarray*}
\frac{ds}{db^{2}}=\frac{1}{2}(\frac{1}{b^{2}}-b^{2}t_{2})s+\frac{1}{2}t_{2}s^{3}.
\end{eqnarray*}
This is a Bernoulli equation which can be rewritten as
\begin{eqnarray*}
\frac{d}{db^{2}}\left(\frac{1}{s^{2}}\right)=\left(b^{2}t_{2}-\frac{1}{b^{2}}\right)\frac{1}{s^{2}}-t_{2}.
\end{eqnarray*}
This is a linear 1-order ODE of $\frac{1}{s^{2}}$.
One can easily get its solution
\begin{eqnarray}\label{chara2}
\frac{1}{s^{2}}=e^{\int(b^{2}t_{2}-\frac{1}{b^{2}})db^{2}}\left[\tilde{c}_{1}-\int t_{2}e^{\int(\frac{1}{b^{2}}-b^{2}t_{2})db^{2}}db^{2}\right],
\end{eqnarray}
where $\tilde{c}_{1}$ is a constant. Hence, by (\ref{chara2}), one independent integral of Eq. (\ref{chara}) is
\begin{eqnarray}\label{chara3}
\frac{s^{2}}{e^{\int(\frac{1}{b^{2}}-b^{2}t_{2}))db^{2}}+s^{2}\int t_{2}e^{\int(\frac{1}{b^{2}}-b^{2}t_{2})db^{2}}db^{2}}=\frac{1}{\tilde{c}_{1}}.
\end{eqnarray}
Note that the characteristic equation (\ref{chara}) is equivalent to
\begin{eqnarray}\label{chara4}
\frac{db^{2}}{1}=\frac{d\ln s}{\frac{1}{2}[\frac{1}{b^{2}}-(b^{2}-s^{2})t_{2}]}=\frac{d\ln\phi}{\frac{1}{2}(-\frac{1}{b^{2}}+t_{2}s^{2})}.
\end{eqnarray}
Eq. (\ref{chara4}) implies
\begin{eqnarray}\label{chara5}
\frac{db^{2}}{1}=\frac{d\ln s-d\ln\phi}{\frac{1}{b^{2}}-\frac{1}{2}b^{2}t_{2}}
\end{eqnarray}
By integrating Eq. (\ref{chara5}), we obtain another independent integral of Eq. (\ref{chara})
\begin{eqnarray}\label{chara6}
\ln \frac{s}{\phi}-\int(\frac{1}{b^{2}}-\frac{1}{2}b^{2}t_{2})db^{2}=\tilde{c}_{2},
\end{eqnarray}
where $\tilde{c}_{2}$ is a constant. Hence, the general solution of Eq. (\ref{1order}) is
\begin{eqnarray}\label{generalsolution}
\Phi\left(\frac{s^{2}}{e^{\int(\frac{1}{b^{2}}-b^{2}t_{2}))db^{2}}+s^{2}\int t_{2}e^{\int(\frac{1}{b^{2}}-b^{2}t_{2})db^{2}}db^{2}}, \ln \frac{s}{\phi}-\int(\frac{1}{b^{2}}-\frac{1}{2}b^{2}t_{2})db^{2}\right)=0,
\end{eqnarray}
where $\Phi(\xi,\eta)$ is any continuously differentiable function.
Suppose $\Phi'_{\eta}\neq0$, then we can solve from (\ref{generalsolution}) that
 \begin{eqnarray}\label{generalsolution1}
 \phi=\varphi(\frac{s^{2}}{e^{\int(\frac{1}{b^{2}}-b^{2}t_{2}))db^{2}}+s^{2}\int t_{2}e^{\int(\frac{1}{b^{2}}-b^{2}t_{2})db^{2}}db^{2}})e^{\int(\frac{1}{2}b^{2}t_{2}-\frac{1}{b^{2}})db^{2}}s,
 \end{eqnarray}
 where $\varphi(\cdot)$ is any positive continuously differentiable function. Hence, the corresponding general $\ab$-metric is
 \begin{eqnarray*}
 F=\alpha \varphi(\frac{s^{2}}{e^{\int(\frac{1}{b^{2}}-b^{2}t_{2}))db^{2}}+s^{2}\int t_{2}e^{\int(\frac{1}{b^{2}}-b^{2}t_{2})db^{2}}db^{2}})e^{\int(\frac{1}{2}b^{2}t_{2}-\frac{1}{b^{2}})db^{2}}s.
 \end{eqnarray*}

\end{proof}

\noindent Hongmei Zhu\\
College of Mathematics and Information Science, Henan Normal University, Xinxiang, 453007, P.R. China\\
zhm403@163.com
\end{document}